\documentstyle[12pt]{article}
\def\ptl{\partial}

\def\be{\beta}
\def\G{\Gamma}

\def\vs{\varsigma}

\def\ad{\mbox{ad}}

\def\ol{\overline}
\def\ul{\underline}

\def\for{\mbox{for}}

\def\rar{\rightarrow}
\def\bs{\backslash}
\def\der{\mbox{Der}\:}

\def\a{\alpha}
\def\b{\beta}
\def\N{\hbox{$I\hskip -4pt N$}}
\def\F{\hbox{$I\hskip -4pt F$}}
\def\Z{\hbox{$Z\hskip -5.2pt Z$}}
\def\sZ{\hbox{$\sc Z\hskip -4.2pt Z$}}

\def\qed{\hfill \hfill \ifhmode\unskip\nobreak\fi\ifmmode\ifinner
\else\hskip5pt\fi\fi\hbox{\hskip5pt\vrule width4pt height6pt
depth1.5pt\hskip 1 pt}}
\def\qed{\hfill$\Box$}

\def\g{\gamma}
\def\G{\Gamma}

\def\si{\sigma}
\def\sc{\scriptstyle}
\def\ssc{\scriptscriptstyle}
\def\dis{\displaystyle}
\def\cl{\centerline}
\def\DD{\hbox{$I\hskip -4pt D$}}

\def\sF{\hbox{$\sc I\hskip -2.5pt F$}}

\def\ol{\overline}
\def\ul{\underline}

\def\rar{\rightarrow}
\def\Rar{\Rightarrow}

\def\Lar{\Leftarrow}

\def\Rla{\Leftrightarrow}
\def\bs{\backslash}
\def\hs{\hspace*}
\def\vs{\vspace*}
\def\rb{\raisebox}

\def\ni{\noindent}
\def\hi{\hangindent}
\def\ha{\hangafter}
\def\mathbb#1{{\mbox{c$\!\!\!\!$}\cal #1}}

\def\AA{{\cal A}}
\def\DD{{\cal D}}
\def\BB{{\cal B}}
\def\ii{{\bf i}}
\def\jj{{\bf j}}
\def\kk{{\bf k}}
\def\mm{{\bf m}}
\def\sone{{1\hskip -5.3pt 1}}
\def\one{{1\hskip -4.7pt 1}}
 
\begin{document}
\par\
\par\
\par\ni
\cl{\bf Structure of the Lie algebras related to those of block}
\par\cl{(Appeared in {\it Comm.~Alg.} {\bf30} (2002), 3205-3226)}
\vs{5pt} {\small\par \cl{\bf Yucai Su} \cl{Department of
Mathematics, Shanghai Jiaotong University, Shanghai 200030, China}
\par
\centerline{and} 
\par \cl{\bf Jianhua Zhou} \cl{Department
of Applied Mathematics, Southeast University, Nanjing 210096,
China}
\par\
\vs{-5pt}\par \cl{\bf ABSTRACT}
\par
We determine the isomorphism classes of the first family of
infinite dimensional simple Lie algebras recently
introduced by Xu. The structure space of these algebras
is given explicitly. The derivations of these algebras are also determined.
\par\ni
\par\ni
 }\par\ \vs{-5pt}\par \cl{\bf 1. \ INTRODUCTION}
\par
Block~(1) introduced a class of infinite dimensional simple Lie algebras
over a field of characteristic zero. Recently there have appeared many
papers (see, for example, Refs.~(2)-(4),~(6),~(10)-(14)) on infinite
dimensional simple Lie algebras which can be viewed as generalizations
of Block algebras. These Lie algebras are constructed from pairs
$(\AA,\DD)$ consisting of a commutative associative algebra $\AA$
with an identity
element and a finite dimensional Abelian derivation subalgebra $\DD$
such that $\AA$ is $\DD$-simple. Such pairs $(\AA,\DD)$ with $\DD$ being
locally finite were classified in Ref.~(8). Based on this classification
it was very natural for Xu~(12) to give generalizations of Block algebras.
\par
In this paper, we shall determine the derivation algebras and the
isomorphism classes of the first family of
infinite dimensional simple Lie algebras introduced by Xu~(12).
Below we shall introduce the normalized form of these algebras.
Let $\F$ be a field of characteristic 0. Denote by $\Z$ the ring of
integers and by $\N$ the non-negative integers $\{0,1,2,...\}$.
\par
Pick $J_p\in\{\{0\},\N\}$ for $p=1,2$,
and set $J=J_1\times J_2$.
Denote by $\pi_p$ the projection from $\F^2$ to the $p$th coordinate,
that is,
$\pi_p(\a)=\a_p$ for $\a=(\a_1,\a_2)\in\F^2,p=1,2.$
Take an additive subgroup $\G$ of $\F^2$ such that
$$
\pi_p(\G)\ne\{0\}\mbox{ \ if \ }J_p=\{0\}\ \ \for\ \ p=1,2.
\eqno(1.1)$$
For $a\in\F$, we denote
$$
a_{[1]}=(a,0),\ a_{[2]}=(0,a)\in\F^2.
\eqno(1.2)$$
Let $\AA_2=\AA_2(\G,J)$ be the semigroup algebra with a basis
$\{x^{\a,\ii}\,|\,(\a,\ii)\in\G\times J\}$
and the algebraic operation $\cdot$ defined by:
$$
x^{\a,\ii}\cdot x^{\b,\jj}=x^{\a+\b,\ii+\jj}\ \ \for\ \
(\a,\ii),(\b,\jj)\in\G\times J.
\eqno(1.3)$$
Then $\AA_2$ forms a commutative associative algebra with an identity
element $1=x^{0,0}$. We define the
derivations $\ptl_p$ of $\AA_2$ by
$$
\ptl_p(x^{\a,\ii})=\a_p x^{\a,\ii}+i_p x^{\a,\ii-1_{[p]}}\ \ \for\ \
(\a,\ii)\in\G\times J,\ p=1,2,
\eqno(1.4)$$
where we adopt the convention that if a notion is not defined but technically
appears in an expression, we always treat it as zero; for instance,
$x^{\a,-1_{[1]}}=0$ for any $\a\in\G$.
\par
We define the following algebraic operation $[\cdot,\cdot]$
on $\AA_2=\AA_2(\G,J)$:
$$
[u,v]=\ptl_1(u)\ptl_2(v)-\ptl_1(v)\ptl_2(u)+u\ptl_1(v)-v\ptl_1(u)
\ \for\ u,v\in\AA_2.
\eqno(1.5)$$
Denote
$$
\si_1=(0,1),\ \ \si_2=(0,2).
\eqno(1.6)$$
We treat $x^{\si_i,0}=0$ if $\si_i\notin\G$ for $i=1,2.$
Then $x^{\si_1,0}$ is a central element of $\AA_2$. Form a quotient Lie
algebra
$$
\BB_2=\BB_2(\G,J)=\AA_2/\F x^{\si_1,0},
\eqno(1.7)$$
whose induced Lie bracket is still denoted by $[\cdot,\cdot]$.
\par\ni
{\bf Theorem 1.1}. {\it The Lie algebra $(\BB_2,[\cdot,\cdot])$ is simple
if $J\ne\{0\}$ or $\si_2\notin\G$. If $J=\{0\}$ and $\si_2\in\G$, then
$\BB_2^{(1)}=[\BB_2,\BB_2]$ is simple and $\BB_2=\BB_2^{(1)}
\oplus(\F x^{\si_1,0}+\F x^{\si_2,0})$.}
\qed\par
The above theorem was due to Block~(1) when $J=\{0\}$ and $\si_2\notin\G$,
due to Dokovic and Zhao~(2) when $J=\{0\}$ and $\si_2\in\G$, and due to
Xu~(12) when $J\ne\{0\}$.
\par
We shall simply denote the simple Lie algebra mentioned
in Theorem 1.1 by $\BB=\BB(\G,J)$. The Lie algebras $\BB$ are the
normalized forms of the first family of the simple Lie algebras
introduced by Xu~(12).
We shall use notation $x^{\a,\ii}$ again to denote elements
in $\BB$. In particular, we have $x^{\si_1,0}=0$.
We set $\G^\#=\G$ if $J\ne\{0\}$ and
$\G^\#=\G\bs\{\si_1,\si_2\}$ otherwise, then $\BB$ has a basis
$$
B=\{x^{\a,\ii}\,|\,(\a,\ii)\in\G^\#\times J,(\a,\ii)\ne(\si_1,0)\}.
\eqno(1.8)$$
\hs{3ex}
Denote
by $GL_2$ the group of invertible
$2\times 2$ matrices with entries in $\F$.
Define an action of $GL_2$ on $\F^2$ by
$g(\a)=\a g^{-1}$ for $\a\in\F^2,g\in GL_2$.
For any additive subgroup $\Upsilon$
of $\F^2$ and $g\in GL_2$, the set
$$
g(\Upsilon)=\{g(\a)\mid \a\in \Upsilon\},
\eqno(1.9)$$
also forms an additive subgroup of $\F^2$.
Denote
by $\Omega_1$ and $\Omega_2$ and $\Omega_3$ the sets of subgroups $\G$ of
$\F^2$ satisfying $\pi_1(\G)\!\ne\!\{0\}\!\ne\!\pi_2(\G)$ and
$\pi_2(\G)\!\ne\!\{0\}$ and $\pi_1(\G)\!\ne\!\{0\}$ respectively,
by $\Omega_4$ the sets of subgroups $\G$ of $\F^2$.
Denote by $G_1$ and $G_2$ the groups of invertible $2\!\times\!2$ matrices of
the forms
$\pmatrix{a\!\!\!\!&b\cr0\!\!\!\!&1\cr}$
and
$\pmatrix{a\!\!\!\!&0\cr0\!\!\!\!&1\cr}$
respectively.
We define the moduli space ${\cal M}_1\!=\!\Omega_1/G_1$,
which is the set of $G_1$-orbits in $\Omega_1$ under the action (1.9),
and define ${\cal M}_2\!=\!\Omega_2/G_2$,
${\cal M}_3\!=\!\Omega_3/G_1$,
${\cal M}_4\!=\!\Omega_4/G_1.$
The main theorem of this paper is the following.
\par\ni
{\bf Theorem 1.2}. {\it The Lie algebras
$\BB=\BB(\G,J)$ and $\BB'=\BB(\G',J')$ are isomorphic if and only if
$J=J'$ and there exist $a,b\in\F,a\ne0$ such that
$b=0$ if $J_1\ne\{0\}=J_2$ and
the map
$\phi: (\b_1,\b_2)\mapsto(a\b_1,\b_2+b\b_1)$ is a group isomorphism
$\G\cong\G'$  (in particular, if
$\pi_1(\G)
=\{0\}$, then $\BB\cong\BB'\Rla (\G,J)=(\G',J')$). Thus there is a bijection
between the isomorphism classes of the simple Lie algebras in Theorem 1.1
and the following set:
$$
{\cal M}=\{(i,\omega)\,|\,i=1,2,3,4,\,\omega\in{\cal M}_i\}.
\eqno(1.10)$$ In other word, ${\cal M}$ is the structure space of
these simple Lie algebras. }\par\ni {\bf Remark 1.3}. {\it We can
drop the condition $\pi_2(\G)+J_2\ne\{0\}$, since if
$\pi_2(\G)=J_2=\{0\}$, then the Lie algebra $\BB(\G,J)$ is a known
rank-one simple Lie algebra of generalized Witt type in
Refs.~(5),~(8),~(9). In this case, the structure space is the set
which is the union of two copies of ${\cal M}_3$ and two copies of
${\cal M}_4$.} \qed\par\ \vs{-5pt}\par \cl{\bf 2. \ STRUCTURE OF
THE DERIVATION ALGEBRAS}
\par
In this section,
we shall determine the structure of the derivation algebra
of the Lie algebra $\BB=\BB(\G,J)$.
Recall that a derivation $d$ of the Lie algebra $\BB$ is a linear
transformation on $\BB$ such that
$$
d([u_1,u_2])=[d(u_1),u_2]+[u_1,d(u_2)]\ \ \for\ \ u_1,u_2\in\BB.
\eqno(2.1)$$
Denote by $\der\BB$ the space of the derivations of $\BB$,
which is a Lie algebra. Moreover, $\ad_{\cal B}$ is an ideal.
Elements in $\ad_{\cal B}$ are called
{\it inner derivations}, while elements in $\der\BB\bs\ad_{\cal B}$ are
called {\it outer derivations}.
Note that
$\BB=\oplus_{\a\in\G}\BB_\a$ is a $\G$-graded Lie algebra with
$$
\BB_\a={\rm span}\{x^{\a,{\bf i}}\,|\,\ii\in J\}\ \for\ \a\in\G.
\eqno(2.2)$$
Define a total order on $J$ such that $\ii>\jj$ if $|\ii|>|\jj|$ or
$|\ii|=|\jj|$ but $i_1>j_1$,
where $|\ii|=i_1+i_2$ is called the {\it level} of $\ii$. Denote
$$
\BB_\a^{[\ii]}={\rm span}\{x^{\a,\jj}\,|\,\jj\le \ii\},\
\BB_\a^{(\ii)}={\rm span}\{x^{\a,\jj}\,|\,\jj< \ii\},
\eqno(2.3)$$
for $(\a,\ii)\in\G\times J$.
For $\a\in\G$, we set
$$
(\der\BB)_\a=\{d\in\der\BB\mid d(\BB_\b)
\subset\BB_{\a+\b}\;\for\;\b\in\G\}.
\eqno(2.4)$$
\hs{3ex}
Fix an element $\a\in\G$. Consider a nonzero element
$$
d\in(\der\BB)_\a\mbox{ \ such that \ }
d(\BB^{[\jj]})\subset\BB^{[\ii+\jj]}
\mbox{ \ for any \ }\jj\in J,
\eqno(2.5)$$
where $\ii\in\Z^2$ is a fixed element. Define the order in $\Z^2$ similarly
as in $J$.
We assume
$$
\ii\mbox{ \ in (2.5) is the minimal element satisfying the condition}.
\eqno(2.6)$$
\hs{3ex}
For any $(\b,\jj)\in\G^\#\times J$ with $(\b,\jj)\ne(\si_1,0)$,
we define $e_{\b,\jj}\in\F$ by
$$
d(x^{\b,\jj})\equiv e_{\b,\jj}x^{\a+\b,\ii+\jj}
\ ({\rm mod\,}\BB^{(\ii+\jj)}),
\eqno(2.7)$$
where if $(\a+\b,\ii+\jj)=(\si_1,0)$, we shall assume that $e_{\b,\jj}=0$.
We rewrite (1.5) as
$$
\matrix{
[x^{\a,\ii},x^{\b,\jj}]
=\!\!\!\!&
(\a_1(\b_2\!-\!1)-\b_1(\a_2\!-\!1))x^{\a+\b,\ii+\jj}
\vs{4pt}\hfill\cr&
+(i_1(\b_2\!-\!1)-j_1(\a_2\!-\!1))x^{\a+\b,\ii+\jj-1_{[1]}}
\vs{4pt}\hfill\cr&
+(\a_1j_2-\b_1i_2)x^{\a+\b,\ii+\jj-1_{[2]}}
+(i_1j_2-j_1i_2)x^{\a+\b,\ii+\jj-1_{[1]}-1{[2]}},
\hfill\cr}
\eqno(2.8)$$
for $x^{\a,\ii},x^{\b,\jj}\in B$.
Then we have
$$
[x^{\b,\jj},x^{\g,\kk}]\equiv (\b_1(\g_2-1)-\g_1(\b_2-1))x^{\b+\g,\jj+\kk}
\ ({\rm mod\,}\BB^{(\jj+\kk)}),
\eqno(2.9)$$
for $(\b,\jj),(\g,\kk)\in\G^\#\times J$ with
$(\b,\jj),(\g,\kk),(\b+\g,\ii+\kk)\ne(\si_1,0)$
(we shall always assume that this condition
is satisfied in the following discussion).
\mbox{Applying} $d$ to (2.9), we get
$$
\matrix{
((\a_1+\b_1)(\g_2-1)-\g_1(\a_2+\b_2-1))e_{\b,\jj}
\vs{4pt}\hfill\cr\ \ \ \ \ \
+(\b_1(\a_2+\g_2-1)-(\a_1+\g_1)(\b_2-1))e_{\g,\kk}
\vs{4pt}\hfill\cr
=(\b_1(\g_2-1)-\g_1(\b_2-1))e_{\b+\g,\jj+\kk}.
\hfill\cr}
\eqno(2.10)$$
\par\ni
{\bf Lemma 2.1}. {\it If $\a_1\ne0$, then the derivation $d$ in (2.5) is
an inner derivation of $\BB$.}
\par\ni
{\it Proof}.
Taking $\g=0,\,\kk=0$ in (2.10), we obtain
$-(\a_1+\b_1)e_{\b,\jj}
+(\b_1(\a_2-1)-\a_1(\b_2-1))e_{0,0}
=-\b_1 e_{\b,\jj}$, i.e.,
$$
e_{\b,\jj}=\a_1^{-1}(\b_1(\a_2-1)-\a_1(\b_2-1))e_{0,0}.
\eqno(2.11)$$
This shows that $e_{\b,\jj}$ does not depend on $\jj$ for any $\b$. Letting
$\jj=0$ in (2.7), we obtain that $\ii\in J$ by the assumption (2.6).
Set $u=\a_1^{-1}e_{0,0}x^{\a,\ii}\in\BB$ and let
$d'=d-\ad_u.$
By (2.8) and (2.11), we obtain that
$d'(x^{\b,\jj})\in\BB^{(\ii+\jj)}$ for all
$(\b,\jj)\in\G^\#\times J$. By induction on $\ii$ in (2.6),
$d'$ is an inner derivation. Hence $d$ is an inner derivation.
\qed\par
Note that $\BB(\G,J)$ is a subalgebra of $\BB(\G,\N^2)$. If $\si_1\in\G$
and $J_p=\{0\}$ for some $p=1,2$, then $x^{\si_1,1_{[p]}}\notin\BB$ but
$[x^{\si_1,1_{[p]}},\BB]\subset\BB$.
Thus we can define outer derivations $d_1,\ol d_1$ as follows:
\par\ni\hs{3pt}$
\matrix{
d_1\!=\!\ad_{x^{\si_1,1_{[2]}}}|_\BB:
x^{\b,\jj}\!\mapsto\!-\b_1x^{\si_1+\b,\jj}
\!-\!j_1x^{\si_1+\b,\jj-1_{[1]}}\hfill
\!\!\!\!&\mbox{if }\si_1\!\in\!\G,J_2\!=\!\{0\},
\vs{4pt}\hfill\cr
\ol d_1\!=\!\ad_{x^{\si_1,1_{[1]}}}|_\BB:
x^{\b,\jj}\!\mapsto\!(\b_2-1)x^{\si_1+\b,\jj}
\!+\!j_2x^{\si_1+\b,\jj-1_{[2]}}\hfill
\!\!\!\!&\mbox{if }\si_1\!\in\!\G,J_1\!=\!\{0\}.
\hfill\cr
}
$\hfill(2.12)\par\ni
We have another outer derivation $d_2$ defined by
$$
d_2=\ad_{x^{\si_2,0}}|_\BB:x^{\b,0}\mapsto -\b_1x^{\si_2+\b,0}
\mbox{ \ if \ }\si_2\in\G,\,J=\{0\}.
\eqno(2.12)'$$
We set $d_1,\ol d_1,d_2$ to be zero whenever they are not defined.
\par\ni
{\bf Lemma 2.2}. {\it If $\a_1=0$ and $\a_2\ne0$,
then the derivation $d$ in (2.5) has the form $d=\ad_u+c_1d+c'_1\ol d_1+c_2d_2$
for $u\in\BB,c_1,c'_1,c_2\in\F$.}
\par\ni
{\it Proof}.
Now (2.10) take the following form
$$
\matrix{
(\b_1(\g_2{\sc\!}-{\sc\!}1)
-\g_1(\a_2{\sc\!}+{\sc\!}\b_2{\sc\!}-{\sc\!}1))e_{\b,\jj}
+(\b_1(\a_2{\sc\!}+{\sc\!}\g_2{\sc\!}-{\sc\!}1)
-\g_1(\b_2{\sc\!}-{\sc\!}1))e_{\g,\kk}
\vs{4pt}\hfill\cr
=(\b_1(\g_2-1)-\g_1(\b_2-1))e_{\b+\g,\jj+\kk}.
\hfill\cr}
\eqno(2.10)'$$
We consider the following cases.
\par
{\it Case 1}. $\pi_1(\G)\ne\{0\}$ and $\a_2\ne1$.
\par
In the following, we take $\b\in\G$ to be arbitrary with $\b_1\ne0$.
First by letting $\g=\kk=0$ in $(2.10)'$, we obtain $e_{0,0}=0$.
Now take $\g\in\G$ to be arbitrary such that $\g_1=0,\g_2\ne0,1$
(such $\g$ must exist, for example, we can take $\g$ to be a
multiple of $\a$).
Letting $\kk=0$, then in $(2.10)'$, by canceling the common factor
$\b_1$, we obtain
$$
(\g_2-1)(e_{\b,\jj}-e_{\b+\g,\jj})+(\a_2+\g_2-1)e_{\g,0}=0.
\eqno(2.13)$$
Substituting $\b$ by $\b+\g$ and $\g$ by $-\g$, we obtain
$(\g_2+1)(e_{\b,\jj}-e_{\b+\g,\jj})+(\a_2-\g_2-1)e_{-\g,0}=0.$
From this and (2.13), we obtain
$$
(\g_2+1)(\a_2+\g_2-1)e_{\g,0}=(\g_2-1)(\a_2-\g_2-1)e_{-\g,0}.
\eqno(2.14)$$
Substituting $\g$ by $-\b+\g$ in $(2.10)'$ and letting $\jj=\kk=0$, we obtain
$$
\matrix{
\hfill(\g_2+\a_2-2)(e_{\b,0}+e_{-\b+\g,0})
\!\!\!&=(\g_2-2)e_{\g,0},
\vs{4pt}\hfill\cr
\hfill(-\g_2+\a_2-2)(e_{\b,0}+e_{-\b-\g,0})
\!\!\!&=(-\g_2-2)e_{-\g,0}.
\hfill\cr}
\eqno\matrix{(2.15)\vs{4pt}\cr(2.16)\cr}$$
where (2.16) is obtained from (2.15) by replacing $\g$ by $-\g$.
\par
{\bf Claim 1}. $e_{\g,0}=0$
if $\a_2=2$ or $2\g_2^2+3\a_2-5\ne0$.
\par
First assume that $\a_2=2$. Adding (2.15) to (2.16) and re-denoting
$-\b-\g$ by $\b$ gives
$\g_2(e_{\b,0}-e_{\b+2\g})=(\g_2+2)e_{-\g,0}-(\g_2-2)e_{\g,0}.$
Multiplying it by $(\g_2-1)^2$ and using (2.14) (noting that $\a_2=2$)
to substitute $e_{-\g,0}$, we obtain
$$
\g_2(\g_2-1)^2
(e_{\b,0}-e_{\b+2\g})=-((\g_2+2)(\g_2+1)^2+(\g_2-2)(\g_2-1)^2)e_{\g,0}.
\eqno(2.17)$$
Replacing $\b$ by $\b+\g$ in (2.13) and adding the result to (2.13), we obtain
$(\g_2-1)(e_{\b,\jj}-e_{\b+2\g,\jj})=-2(\a_2+\g_2-1)e_{\g,0}.$ Comparing this
with (2.17), since the determinant
$$
\left|\matrix{\g_2(\g_2-1)^2&\!\!(\g_2{\sc\!}+{\sc\!}2)
(\g_2{\sc\!}+{\sc\!}1)^2{\sc\!}+{\sc\!}(\g_2{\sc\!}-{\sc\!}2)
(\g_2{\sc\!}-{\sc\!}1)^2
\vs{4pt}\cr
(\g_2-1)&\!\!2(\a_2+\g_2-1)
\cr}\right|
=-12\g_2(\g_2-1)\ne0,
\eqno(2.18)$$
we obtain $e_{\g,0}=0$. Now assume that $\a_2\ne2$.
Setting $\g$ to be zero in (2.15) gives $e_{-\b,0}=-e_{\b,0}$.
Thus (2.15) becomes
$(\g_2+\a_2-2)(e_{\b,0}-e_{\b-\g,0})=(\g_2-2)e_{\g,0}.$
Replacing $\g$ by $-\g$ and $\b$ by $\b-\g$ gives
$(\g_2-\a_2+2)(e_{\b,0}-e_{\b-\g,0})=-(\g_2+2)e_{-\g,0}.$
From these two formulas, we obtain
$(\g_2-\a_2+2)(\g_2-2)e_{\g,0}.=-(\g_2+\a_2-2)(\g_2+2)e_{-\g,0}.$
Comparing this with (2.14), since
$$
\left|\matrix{(\g_2+1)(\a_2+\g_2-1)&(\g_2-1)(\a_2-\g_2-1)
\vs{4pt}\cr
(\g_2-\a_2+2)(\g_2-2)&-(\g_2+\a_2-2)(\g_2+2)
\cr}\right|
=-2\a_2\g_2(2\g_2^2+3\a_2-5)\ne0,
$$
we obtain $e_{\g,0}=0$. This proves Claim 1.
\par
Claim 1 and (2.13) show that $e_{\b,\ii}=e_{\b+\g,\ii}$ under the
condition in Claim 1, and from this, we can easily deduce that
$e_{\b,\ii}=e_{\b+\g,\ii}$ holds for all $\g$ with
$\g_1=0,\g_2\ne0,1$ (and thus condition in Claim 1 can be removed). Then
by substituting $\g$ by $\b+\g$ in $(2.10)'$, we obtain that
$e_{2\b,2\ii}=2e_{\b,\ii}$.
Induction on $m\in\N$ gives $e_{m\b,m\ii}=me_{\b,\ii}$. Let
$\b'\in\G$ with $\b'_1\ne0$. In $(2.10)'$, by replacing
$\b,\ii,\g,\kk$ by $m\b,m\ii,m\b',m\kk$ respectively,
we obtain
$$
\matrix{
(\b_1(m\b'_2-1)-\b'_1(\a_2+m\b_2-1))e_{\b,\jj}
\vs{4pt}\hfill\cr\ \ \ \ \ \
+(\b_1(\a_2+m\b'_2-1)-\b'_1(m\b_2-1))e_{\b',\kk}
\vs{4pt}\hfill\cr
=(\b_1(m\b'_2-1)-\b'_1(m\b_2-1))e_{\b+\b',\jj+\kk}.
\hfill\cr}
\eqno(2.19)$$
Computing the coefficient of $m$ implies $e_{\b+\b',\jj+\kk}
=e_{\b,\jj}+e_{\b',\kk}$. Using this in (2.19) gives
$\b'_1e_{\b,\jj}=\b_1e_{\b',\kk}.$
Fixing $\b'$ and $\kk$
shows that $e_{\b,\jj}=f\b_1$ for some $f\in\F$ and for all $(\b,\jj)
\in\G^\#\times J$.
Thus as in the proof of Lemma 2.1, $\ii\in J$, and by setting
$d'=d-f\ad_{x^{\a,\ii}}|_{\BB}$
(note that if $J=\{0\}$ and $\a=\si_2$, then $\ad_{x^{\a,\ii}}|_{\BB}
=d_2$ is not inner)
and induction on $\ii$, we see that $d$ has the
required form.
\par
{\it Case 2}. $\pi_1(\G)\ne\{0\}$ and $\a_2=1$.
\par
Then $(2.10)'$ becomes
$$
\matrix{
(\b_1(\g_2-1)-\g_1\b_2)e_{\b,\jj}
+(\b_1\g_2-\g_1(\b_2-1))e_{\g,\kk}
\vs{4pt}\hfill\cr
=(\b_1(\g_2-1)-\g_1(\b_2-1))e_{\b+\g,\jj+\kk}.
\hfill\cr}
\eqno(2.10)''$$
Setting $\g=0$, we see that $e_{\b,\ii}$ does not depend on $\ii$ if
$\b_1\ne0$, and then $(2.10)''$ also shows that $e_{\b,\ii}$ does not
depend on $\ii$ if $\b_1=0$. Thus $\ii\in J$. Denote $e_\b=e_{\b,\ii}$.
\par
First assume that $e_0\ne0$. If $i_2\ne0$, taking
$\b\in\G$ with $\b_1\ne0$,
we have
$$
\matrix{
-i_2\b_1e_0x^{\b+\si_1,\ii-1_{[2]}}
\!\!\!\!&
\equiv[e_0x^{\si_1,\ii},x^{\b,0}]
\equiv[d(1),x^{\b,0}]
\vs{4pt}\hfill\cr&
=d([1,x^{\b,0}])-[1,d(x^{\b,0})]
\vs{4pt}\hfill\cr&
\equiv(\b_1-\ptl_1)d(x^{\b,0})
\equiv0\,(\,\mbox{mod}\,\BB^{[\ii-1_{[1]}]}\,),
\hfill\cr}
\eqno(2.20)$$
a contradiction. Thus $i_2=0$.
Let $d'=\ad_{x^{\si_1,\ii+1_{[1]}}}|_\BB$.
Note that $d'$ is either an inner derivation if $J_1\ne\{0\}$ or else
$d'=\ol d_1$ (cf.~(2.12)), and $d'(\BB^{[\jj]})\subset\BB^{[\ii+\jj]}$. Thus
by replacing $d$ by $d-fd'$ for some $f\in\F$, we can suppose $e_0=0$.
Now exactly as in Case 1, we have $e_\b=f\b_1$ for some $f\in\F$ and
by setting $d'=d-fd''$ and induction on $\ii$,
we see that $d$ has the required form, where $d''=
\ad_{x^{\si_1,\ii+1_{[2]}}}|_\BB$ is either inner or equal to $d_1$.
\par
{\it Case 3}. $\pi_1(\G)=\{0\}$.
\vs{-1pt}\par
Then by (1.1), $J_1\ne\{0\}$. By (2.8), we have
$$
[x^{\b,\jj},x^{\g,\kk}]\equiv
(j_1(\g_2-1)-k_1(\b_2-1))x^{\b+\g,\jj+\kk-1_{[1]}}
\ ({\rm mod\,}\BB^{(\jj+\kk-1_{[1]})}),
\eqno(2.21)$$
for $(\b,\jj),(\g,\kk)\in\G^\#\times J$.
Applying $d$ to (2.21), we get
$$
\matrix{
((i_1+j_1)(\g_2\!-\!1)-k_1(\a_2+\b_2\!-\!1))e_{\b,\jj}
\vs{4pt}\hfill\cr\ \ \ \ \
+(j_1(\a_2+\g_2\!-\!1)-(i_1+k_1)(\b_2\!-\!1))e_{\g,\kk}
\vs{4pt}\hfill\cr
=(j_1(\g_2-1)-k_1(\b_2-1))e_{\b+\g,\jj+\kk-1_{[1]}}.
\hfill\cr}
\eqno(2.22)$$
Letting $\g=0,\kk=0$ gives
$$
-(i_1+j_1)e_{\b,\jj}
+(j_1(\a_2-1)-i_1(\b_2-1))e_{0,0}
=-j_1e_{\b,\jj-1_{[1]}}.
\eqno(2.23)$$
Assume first that $i_1\ne0$.
Since $\BB$ is generated by $\{x^{\b,\jj}\,|\,j_1\le2\}$, we must have
$i_1\ge-2$.
Assume that $i_1=-2$. Then by (2.7),
$e_{\b,\jj}=0$ if $j_1\le1$.
In (2.22), let $\g=0,\kk=1_{[1]}$, we get $\a_2=2$. Let $\jj=1_{[1]},
k_1=2$ in (2.22), we get $(\g_2+1)e_{\g,\kk}=(\g_2-2\b_2+1)e_{\b+\g,\kk}$.
Using this in (2.22),
we can easily deduce $e_{\b,\kk}=0$ if $k_1\le2$. Thus
we obtain $d=0$.
Assume that $i_1>0$. By setting $j_1=0,1,2,...$ in (2.23), we obtain
that
$$
e_{\b,\jj}=-((i_1+1)(\b_2-1)-j_1(\a_2-1))(i_1+1)^{-1}e_{0,0},
\eqno(2.24)$$
for $(\b,\jj)\in\G^\#\times J$. Thus (2.24) shows that $e_{0,0}\ne0$
by (2.6). Then setting
$\jj=0$, by (2.7), we see $\ii\in J$. Thus we can take
$u=-(i_1+1)^{-1}e_{0,0}{\ssc\,}x^{\a,\ii+1_{[1]}}$ and
let $d'=d-\ad_u$ and use induction on $\ii$.
If $i_1=-1$, then as above, we obtain $e_{\b,0}=0$ and (2.22), (2.23) give
$e_{\b,\jj}=j_1f$ for some $f\in\F$ and so $d=f{\ssc\,}\ad_{x^{\a,0}}$.
Finally assume that $i_1=0$.
Letting $\g=0,\kk=1_{[1]}$ in (2.22) gives
$e_{\b,\jj}=((\b_2-1)-j_1(\a_2-1))\a_2^{-1}e_{0,1_{[1]}}.$
Thus the arguments after (2.24) can be used to complete the proof.
\qed\par
Now we shall describe homogeneous derivations of degree 0.
Denote by ${\rm Hom}_{\sZ}(\G,\F)$ the set of additive group homomorphisms
from $\G$ to $\F$. For any $\mu\in{\rm Hom}_{\sZ}(\G,\F)$, we can define
a homogeneous derivation $d_\mu$ of degree 0 by:
$$
d_\mu(x^{\b,\jj})=\mu(\b)x^{\b,\jj}\ \for\ (\b,\jj)\in\G^\#\times J.
\eqno(2.25)$$
In this way, we regard ${\rm Hom}_{\sZ}(\G,\F)$ as a subspace of $\der\BB$.
If $J_2\ne\{0\}$, one can verify that $\ptl_{t_2}: x^{\b,\jj}\mapsto
j_2x^{\b,\jj-1_{[2]}}$ is a derivation. If $J_2=\{0\}$, we set $\ptl_{t_2}=0$.
\par\ni
{\bf Lemma 2.3}. {\it Let $d\ne0$ be a derivation satisfying (2.5) with
$\a=0$. Then $d=ad_u+d_\mu+f\ptl_{t_2}$ for some
$u\in\BB_0,\mu\in\mbox{\it Hom}_{\sZ}(\G,\F),f\in\F$.}
\par\ni
{\it Proof.}
We consider two cases.
\par
{\it Case 1}. $\pi_1(\G)\ne\{0\}$.
\par
In (2.10), since $\a=0$, all coefficients are the same, thus we obtain
$$
e_{\b+\g,\jj+\kk}=e_{\b,\jj}+e_{\g,\kk}
\mbox{ \ if \ }\b_1(\g_2-1)-\g_1(\b_2-1)\ne0.
\eqno(2.26)$$
Take $\b_1\ne0$. By (2.26),
$e_{m\b+n\b,\jj+\kk}=e_{m\b,\jj}+e_{n\b,\kk}$ if $m\ne n$. From
this we get
$$
e_{m\b,m\jj}=me_{\b,\jj}\ \for\ m\in\N\mbox{ if }\b_1\ne0.
\eqno(2.27)$$
From this and (2.26), one obtains that (2.26) holds for all $\b,\g$.
\par
First assume that $i_1\ne0$. Then $J_1\ne\{0\}$.
Let $\mm$ be the {\it leading degree} of $d(1)$, that is,
$d(1)$ can be written as $c_0x^{0,\mm}+$ a linear combination
of $x^{0,\kk}$ with $\kk<\mm$, where $c_0\ne0$,
i.e., $d(1)\in\BB^{[\mm]}_0\bs\BB^{(\mm)}_0$.
If $\mm>\ii-1_{[1]}$, then by (2.8),
we see that the leading degree of $[d(1),x^{\b,0}]$
is ${\bf m}$ if we take $\b$ with $\b_1\ne0$,
but $d([1,x^{\b,0}])-[1,d(x^{\b,0})]=(\b_1-\ptl_1)d(x^{\b,0})
\in\BB^{[\ii-1_{[1]}]}$, a contradiction. Thus ${\bf m}\le\ii-1_{[1]}$
and we can suppose
$d(1)-c_1x^{0,\ii-1_{[1]}}\in\BB^{(\ii-1_{[1]})}$
for some $c_1\in\F$. Now applying $d$ to $[1,x^{\b,\jj}]=
\b_1x^{\b,\jj}+j_1x^{\b,\jj-1_{[1]}}$ and computing the coefficient
of $x^{\b,\ii+\jj-1_{[1]}}$, we obtain
$c_1\b_1+(i_1+j_1)e_{\b,\jj}=j_1e_{\b,\jj-1_{[1]}}.$
Setting $j_1=0,1,2,...$ gives $e_{\b,\jj}=-i_1^{-1}c_1\b_1$. Thus again
$\ii\in J$. Then taking
$u=-i_1^{-1}c_1x^{0,\ii}$ and letting $d'=d-\ad_u$ and by
induction on $\ii$, we obtain the result.
\par
Now assume that $i_1=0$.
Then using the arguments for $i_1\ne0$ gives $e_{0,0}=0$
(note that $e_{0,0}=0\Rla$ the leading degree of $d(1)$ is $<\ii$).
We claim that $e_{\b,\jj}$ does not depend on $j_1$. For this, assume that
$J_1\ne\{0\}$.
For any $\g\in\G^\#$ with $\g_1=0$,
we apply $d$ to $[1,x^{\g,\kk}]=k_1x^{\g,\kk-1_{[1]}}$ and
calculate the coefficient of $x^{\g,\ii+\kk-1_{[1]}}$.
Noting that the leading term of $[x^{0,\jj},x^{\g,\kk}]$ is
$\le\jj+\kk-1_{[1]}$ and that
the leading term of $d(1)$ is $<\ii$, we see that
$[d(1),x^{\g,\kk}]\in\BB^{(\ii+\kk-1_{[1]})}$. Thus we obtain
$k_1e_{\g,\kk}=k_1e_{\g,\kk-1_{[1]}}$. So
$e_{\g,\kk}$ does not depend on $k_1$. In particular, $e_{0,\kk}=0$ if
$k_2=0$. Using this in (2.26) also gives that $e_{\b,\jj}$ does not depend
on $j_1$ for all $\b\in\G^\#$.
If $J_2=\{0\}$, then $\ii=0$ and
$e_\b=e_{\b,\jj}$ does not depend on $\jj$, thus by (2.26),
the map $\mu:\b\mapsto e_\b$ defines an element
$\mu\in{\rm Hom}_{\sZ}(\G,\F)$;
let $d'=d-d_\mu$, then using induction on $\ii$ we see that $d'$, and so,
$d$ has the required form.
Assume that $J_2\ne\{0\}$. By (2.26), we obtain
$$
e_{m\b,\jj}=m e_{\b,0}+j_2e_{0,1_{[2]}}\ \for\ m\in\Z,
(\b,\jj)\in\G^\#\times J.
\eqno(2.28)$$
Let $\b_1\ne0$. Suppose
$$
d(x^{\b,\jj})=e_{\b,\jj}x^{\b,\ii+\jj}
+\sum_{{\bf m}<\ii+\jj}e_{\b,\jj}^{({\bf m})}x^{\b,{\bf m}}.
\eqno(2.29)$$
Then applying $d$ to $[x^{\b,0},x^{\b,1_{[2]}}]=\b_1x^{2\b,0}$
and calculating the coefficient of $x^{2\b,\ii}$,
and noting that
$[x^{\b,\jj},x^{\b,k_{[2]}}]\equiv(k-j_2)\b_1x^{2\b,\jj+(k-1)_{[2]}}
\,({\rm mod\,}\BB^{(\jj+(k-1)_{[2]})}),$
we obtain $(1-i_2)e_{\b,0}+(i_2+1)e_{\b,1_{[2]}}=e_{2\b,0}$.
Using (2.28), this gives $e_{0,1_{[2]}}=0$ if $i_2\ne-1$.
\par
If $i_2=-1$, then (2.7) shows that $e_{\b,0}=0$ for all
$\b\in\G$, so (2.28) means that $d(x^{\b,1_{[2]}})=e_{0,1_{[2]}}x^{\b,0}$,
and one can prove by induction on $\jj$ that
$d(x^{\b,\jj})=j_2e_{0,1_{[2]}}x^{\b,\jj-1_{[2]}}$;
this proves that $d=e_{0,1_{[2]}}\ptl_{t_2}$.
Assume that $i_2\ne-1$. Then $e_{0,1_{[2]}}=0$ and (2.28) shows that
$e_{\b,\jj}$, denoted now by $e_{\b}$, does not depend on $\jj$.
Thus as before, $i_2\ge0$. If $i_2=0$, then by replacing $d$
by $d-d_\mu$ for some $\mu\in{\rm Hom}_{\sZ}(\G,\F)$, we see that $d$
has the required form as before.
Thus suppose $i_2>0$. Suppose the coefficient of
$x^{\b,(i_2-1)_{[2]}}$ in $d(x^{\b,0})$ is $e'_\b$.
Applying $d$ to $[x^{\b,0},x^{\g,0}]=(\b_1(\g_2-1)-\g_1(\b_2-1))x^{\b+\g,0}$,
and noting that
$$
[x^{\b,\jj},x^{\g,0}]\equiv
(\b_1(\g_2{\sc\!}-{\sc\!}1){\sc\!}-{\sc\!}\g_1
(\b_2{\sc\!}-{\sc\!}1))x^{\b+\g,\jj}{\sc\!}-{\sc\!}
\g_1j_2x^{\b+\g,\jj-1_{[2]}}
\,({\rm mod\,}\BB^{(\jj-1_{[2]})}),
\eqno(2.30)$$
by calculating the coefficient of $x^{\b+\g,(i_2-1)_{[2]}}$, we obtain
\par\ni\hs{5pt}$
\matrix{
-i_2\g_1e_\b+(\b_1(\g_2-1)-\g_1(\b_2-1))e'_\b
\vs{4pt}\hfill\cr
+i_2\b_1e_\g
+(\b_1(\g_2-1)-\g_1(\b_2-1))e'_\g
=(\b_1(\g_2{\sc\!}-{\sc\!}1){\sc\!}-{\sc\!}\g_1(\b_2{\sc\!}-{\sc\!}1))e'_{\b+\g}.
\hfill\cr}
$\hfill(2.31)\par\ni
That is
$$
i_2(\b_1e_\g-\g_1e_\b)+
(\b_1(\g_2-1)-\g_1(\b_2-1))(e'_\b+e'_\g-e'_{\b+\g})=0.
\eqno(2.32)$$
Let $\b_1\ne0$. Taking $\g=m\b$, by (2.27), we can deduce
$e'_{m\b}=me'_\b,\,m\in\Z.$
Replacing $\b$ by $\b+\g$ and $\g$ by $-\g$ in (2.32) gives
$$
\matrix{
i_2(-(\b_1+\g_1)e_\g+\g_1(e_\b+e_\g))
\vs{4pt}\hfill\cr
\ \ \ \
+((\b_1+\g_1)(-\g_2-1)+\g_1(\b_2+\g_2-1))(e'_{\b+\g}-e'_\g-e'_\b)=0.
\hfill\cr}
\eqno(2.33)$$
Adding (2.33) to (2.32) gives
$(2\b_1\g_2-(2\b_2-3)\g_1)(e'_{\b+\g}-e'_\g-e'_\b)=0.$
From this and the fact that $e'_{m\b}=me'_\b$,
we deduce that $\b\mapsto e'_\b$ is an additive
function. Using this, (2.32) now gives that $e_\b=f\b_1$ for
some $f\in\F$. Now
applying induction to $d'\!=\!d\!-\!f{\ssc\,}\ad_{x^{0,(i_2)_{[2]}}}$
gives the result. This proves Case 1.
\par
{\it Case 2}. $\pi_1(\G)=\{0\}$ (thus $J_1\ne\{0\}$).
\par
Then we have (2.21), and (2.22) becomes
\par\ni\hs{5pt}$
\matrix{
((i_1{\sc\!}+{\sc\!}j_1)(\g_2{\sc\!}-{\sc\!}1)
-k_1(\b_2{\sc\!}-{\sc\!}1))e_{\b,\jj}
+(j_1(\g_2{\sc\!}-{\sc\!}1)-(i_1{\sc\!}+{\sc\!}k_1)
(\b_2{\sc\!}-{\sc\!}1))e_{\g,\kk}
\vs{4pt}\hfill\cr
=(j_1(\g_2-1)-k_1(\b_2-1))e_{\b+\g,\jj+\kk-1_{[1]}}.
\hfill\cr}
$\hfill(2.34)\par\ni
Letting $\g=\kk=0$ gives
$-(i_1+j_1)e_{\b,\jj}
-(j_1+i_1(\b_2-1))e_{0,0}
=-j_1e_{\b,\jj-1_{[1]}}.$
If $i_1<0$, the arguments are similar to those given after (2.23).
If $i_1>0$, by setting $j_1=0,1,2,...$, we obtain
that $e_{\b,\jj}=-((i_1+1)(\b_2-1)+j_1)(i_1+1)^{-1}e_{0,0}$
for $(\b,\jj)\in\G^\#\times J$, which shows that $e_{0,0}\ne0$ by (2.6),
and then setting
$\jj=0$, by (2.7), we see $\ii\in J$; thus we can take
$u=-(i_1+1)^{-1}e_{0,0}x^{0,\ii+1_{[1]}}$ and
let $d'=d-\ad_u$ and use induction on $\ii$.
Assume that $i_1=0$. Then setting $\g=\kk=0$ in (2.34) gives
$$
e_{\b,\jj}=e_{\b,1_{[1]}+(j_2)_{[2]}}-(j_1-1)e_{0,0}=
e_{\b,(j_2)_{[2]}}-j_1e_{0,0}.
\eqno(2.35)$$
Using this in (2.34) with $i_1=0$, we obtain
$$
e_{\b,1_{[1]}+j_{[2]}}+
e_{\g,1_{[1]}+k_{[2]}}=e_{\b+\g,1_{[1]}+(j+k)_{[2]}}\ \for\ j,k\in\N.
\eqno(2.36)$$
If $J_2=\{0\}$, then $\ii=0$,
and the map $\mu:\b\mapsto e_{\b,1_{[1]}}$ defines an element
$\mu\in{\rm Hom}_{\sZ}(\G,\F)$ by (2.36); replacing $d$ by $d-d_{\mu'}$,
where $\mu'=\mu-e_{0,0}\pi_2\in{\rm Hom}_{\sZ}(\G,\F)$ (recall that $\pi_2$
is the projection $\a\mapsto\a_2$),
we can suppose $e_{\b,1_{[1]}}=\b_2e_{0,0}$;
then (2.35) gives $e_{\b,\jj}=(\b_2-1+j_1)e_{0,0}$ and thus by letting
$d'=d-e_{0,0}{\ssc\,}\ad_{x^{1_{[1]}}}$ and by induction on $\ii$ we see that
$d'$, and so, $d$ has the required form.
Assume that $J_2\ne\{0\}$. By (2.35), (2.36), we obtain
\par\ni\hs{5pt}$
\matrix{
e_{m\b,\jj}
\!\!\!\!&=e_{m\b,1_{[1]}+(j_2)_{[2]}}-(j_1{\sc\!}-{\sc\!}1)e_{0,0}
=e_{m\b,1_{[1]}}+e_{0,1_{[1]}+(j_2)_{[2]}}-(j_1{\sc\!}-{\sc\!}1)e_{0,0}
\vs{4pt}\hfill\cr&
=m e_{\b,1_{[1]}}+j_2e_{0,{\sone }}-(j_1-1)e_{0,0}\ \ \ \for\ \ \
m\in\N, (\b,\jj)\in\G^\#\times J,
\hfill\cr}
$\hfill(2.37)\par\ni
where ${\one }=(1,1)$. Write
$d(x^{\b,\jj})$ as in (2.29).
We shall apply $d$ to $[x^{\b,1_{[1]}},x^{\b,{\sone }}]=x^{2\b,1_{[1]}}$
and calculate the coefficient of $x^{2\b,\ii+1_{[1]}}$.
First note that
$$
[x^{\b,\jj},x^{\b,\kk}]=
(\b_2-1)(j_1-k_1)x^{2\b,\jj+\kk-1_{[1]}}+
(j_1k_2-k_1j_2)x^{2\b,\jj+\kk-{\sone }}.
\eqno(2.38)$$
We see that if $x^{2\b,\ii+1_{[1]}}$ appears in the term
$[e_{\b,\jj}^{({\bf m})}
x^{\b,{\bf m}},{\ssc\!}x^{\b,{\sone }}],$ which is a term
in $[d(x^{\b,1_{[1]}}),{\ssc\!}x^{\b,{\sone }}]$, then by (2.38),
either ${\bf m}+{\one }-1_{[1]}=\ii+1_{[1]}$,
or ${\bf m}+{\one }-{\one }=\ii+1_{[1]}$. In the first case,
$m_1=1$ and the coefficient is zero (cf.~the first term of the right-hand
side of (2.38)). In latter case,
${\bf m}=\ii+1_{[1]}$ and the coefficient is
$(1-i_2)e_{\b,1_{[1]}}$ (cf.~the second term of the right-hand
side of (2.38)). Similarly, the coefficient of
$x^{2\b,\ii+1_{[1]}}$ appears in $[x^{\b,1_{[1]}},d(x^{\b,{\sone }})]$
is $(1+i_2)e_{\b,{\sone }}$. Thus
we obtain that $(1-i_2)e_{\b,1_{[1]}}+(1+i_2)e_{\b,{\sone }}=e_{2\b,1_{[1]}}$.
Using (2.37), this gives $e_{0,\sone}=0$ if $i_2\ne-1$.
\par
If $i_2=-1$, (2.7) shows that $e_{\b,0}=0$ for all
$\b\in\G$, and so (2.37) means that $e_{\b,\jj}=j_2e_{0,\sone}$;
then one can prove by induction on $\jj$ that $d=e_{0,\sone}\ptl_{t_2}$.
Assume that $i_2\ne-1$. Then $e_{0,\sone}=0$ and (2.37) gives
$e_{\b,\jj}=e_{\b,1_{[1]}}-(j_1-1)e_{0,0}$.
As before, $i_2$ must be $\ge0$.
If $i_2=0$, then as the arguments after (2.36) by replacing $d$
by $d-d_\mu-f{\ssc\,}\ad_{x^{1_{[1]}}}$ for some $\mu\in{\rm Hom}_{\sZ}(\G,\F),
f\in\F$, we see that $d$ has the required form as before.
Thus suppose $i_2>0$. Suppose the coefficient of
$x^{\b,j_{[1]}+(i_2-1)_{[2]}}$ in $d(x^{\b,j_{[1]}})$ is
$e'_{\b,j}$.
Noting that
$$
[x^{\b,\jj},x^{\g,k_{[1]}}]=
(j_1(\g_2-1)-k(\b_2-1))x^{\b+\g,\jj+(k-1)_{[1]}}-
kj_2x^{\b+\g,\jj+(k-1)_{[1]}-1_{[2]}},
$$
by applying $d$ to $[x^{\b,j_{[1]}},x^{\g,k_{[1]}}]=
(j(\g_2-1)-k(\b_2-1))x^{\b+\g,(j+k-1)_{[1]}}$,
and
computing the coefficient of $x^{\b+\g,(j+k-1)_{[1]}+(i_2-1)_{[2]}}$,
we obtain
$$
\matrix{
-i_2ke_{\b,j_{[1]}}{\sc\!}+{\sc\!}(j(\g_2{\sc\!}-{\sc\!}1)
{\sc\!}-{\sc\!}k(\b_2{\sc\!}-{\sc\!}1))e'_{\b,j}
{\sc\!}+{\sc\!}i_2je_{\g,k_{[1]}}{\sc\!}+{\sc\!}(j(\g_2{\sc\!}-{\sc\!}1)
{\sc\!}-{\sc\!}k(\b_2{\sc\!}-{\sc\!}1))e'_{\g,k}
\vs{4pt}\hfill\cr
=(j(\g_2-1)-k(\b_2-1))e'_{\b+\g,j+k-1}.
\hfill\cr}
$$
That is
$$
i_2(je_{\g,k_{[1]}}-ke_{\b,j_{[1]}})
+(j(\g_2-1)-k(\b_2-1))(e'_{\b,j}+e'_{\g,k}-e'_{\b+\g,j+k-1})=0.
\eqno(2.39)$$
Taking $\g=0,k=1$ gives
$i_2(je_{0,1_{[1]}}-e_{\b,j_{[1]}})
+(-j-(\b_2-1))e'_{0,1}=0\mbox{ for all }j\in\N.$
By computing the coefficients
of $j$ and $j^0$,
using $e_{0,1_{[1]}}=0$ (by setting $\b=0,\jj=1_{[1]}$ in (2.35))
and $e_{\b,j_{[1]}}=e_{\b,0}-je_{0,0}$,
we obtain that
$i_2e_{0,0}-e'_{0,1}=0,\,-i_2e_{\b,0}-(\b_2-1)e'_{0,1}=0.$
Thus we have $e_{\b,0}=-(\b_2-1)e_{0,0}$ and finally we obtain that
$e_{\b,\jj}=-(\b_2-1+j_1)e_{0,0}$ by (2.35), (2.37) and that
$e_{0,\sone}=0$. Thus if we let $u=-e_{0,0}x^{0,1_{[1]}+\ii}$
and set $d'=d-\ad_u$, we obtain the result by induction on $\ii$.
This completes the proof of the lemma.
\qed\par
Note that if $\pi_1(\G)\ne\{0\}$ and $J_1=\{0\}$, then
$\ad_1=d_{\pi_1}\in{\rm Hom}_{\sZ}(\G,\F)$. In this case,
we choose a subspace
${\rm Hom}^*_{\sZ}(\G,\F)$ of ${\rm Hom}_{\sZ}(\G,\F)$ such that
${\rm Hom}_{\sZ}(\G,\F)$
$=\F ad_1\oplus{\rm Hom}^*_{\sZ}(\G,\F)$ as vector
spaces. Otherwise we set ${\rm Hom}^*_{\sZ}(\G,\F)={\rm Hom}_{\sZ}(\G,\F)$.
\par\ni
{\bf Lemma 2.4}. {\it Every homogeneous derivation
$d\in(\der\BB)_\a$ satisfies condition (2.5).}
\par\ni
{\it Proof}.
Let $\G'_0$ be a maximal $\F$-linearly independent subset of $\G$ (which
has at most two elements).
Let $\G'_1=\G'_0\cup(\{\a,\si_1,\si_2\}\cap\G)$.
Let $\G'\subset\G$ be the subgroup of $\G$ generated by $\G'_1$.
Let $\BB'$ be the Lie subalgebra of $\BB$ generated by
$M_0=\{x^{\b,\jj}\,|\,\b\in\G'_1,\jj\in J,\,|\jj|\le4\}.$
Then it is straightforward to check that $\BB'=\BB(\G',J)$,
and $d'=d|_{\BB'}$ is a homogeneous derivation of $\BB'$ of degree $\a$.
Since $M_0$ is a finite set and a derivation is determined by its action
on generators, we see that the derivation $d'$ of $\BB'$ satisfies the
condition in (2.5).
\par
By Lemmas 2.1-3, there exist $u'\in \BB'_\a$ and
$\mu'\in{\rm Hom}^*_{\sZ}(\G',\F),f'_1,...,f'_4\in\F$ such that
$$
d'=d_{u',\mu',f'_1,f'_2,f'_3,f'_4}=
\ad_{u'}+d_{\mu'}+f'_1d_1+f'_2\ol d_1+f'_3d_2+f'_4\ptl_{t_2},
\eqno(2.40)$$
where $(\mu',f'_4)=0$ if $\a\ne0$ and $f'_1=f'_2=0$ if $\a\ne\si_1$ and
$f'_3=0$ if $\a\ne\si_2$.
One can immediately verify that for $u\in\BB,\mu\in{\rm Hom}^*_{\sZ}(\G,\F),
f_1,f_2,f_3,f_4\in\F$, $d_{u,\mu,f_1,f_2,f_3,f_4}|_{\BB(\G_1,J)}=0$ if and
only if $(u,\mu|_{\G_1},f_1,f_2,f_3,f_4)=0$,
where $\G_1$ is any subgroup of $\G$ containing $\G'$.
Hence for any $u\in\BB_\a,\mu\in{\rm Hom}^*_{\sZ}(\G,\F)$,
$d_{u,\mu,f_1,f_2,f_3,f_4}|_{\BB(\G_1,J)}=0$ implies
$(u,\mu|_{\G_1},f_1,f_2,f_3,f_4)=0$.
Let $\G_2$ be the maximal subgroup of $\G$ such that
there exists $u''\in\BB(\G_2,\F),\mu''\in{\rm Hom}^*_{\sZ}(\G_2,\F)$,
$f''_1,f''_2,f''_3,f''_4\in\F$ and
$d|_{\BB(\G_2,J)}=d_{u'',\mu'',f''_1,f''_2,f''_3,f''_4}$
with $\mu''|_{\G'}=\mu'$.
\par
Suppose $\G_2\ne\G$. Take $\b\in\G\bs\G_2$. Let $\G_3$ be the subgroup of
$\G$ generated by $\G'$ and $\b$. Then $\G_3$ is finitely generated. Thus
there exist $u_3,\mu_3,f^{(3)}_1,f^{(3)}_2,f^{(3)}_3$,
$f^{(3)}_4$ such that
$d|_{\BB(\G_3,J)}=d_{u_3,\mu_3,f^{(3)}_1,f^{(3)}_2,f^{(3)}_3,f^{(3)}_4}$.
Then we have
$$
d_{u_3-u',\mu_3-\mu',f^{(3)}_1-f'_1,f^{(3)}_2-f'_2,
f^{(3)}_3-f'_3,f^{(3)}_4-f'_4}|_{\BB(\G',J)}
=d|_{\BB(\G',J)}-d|_{\BB(\G',J)}=0.
\eqno(2.41)$$
Thus $u_3=u',f^{(3)}_i=f'_i,i=1,2,3,4$
and $\mu_3|_{\G'}=\mu'|_{\G'}$. Similarly,
$\mu''|_{\G_2\cap\G_3}=\mu_3|_{\G_2\cap\G_3}$.
Let $\G_4$ be the subgroup of $\G$ generated by $\G_2$ and $\b$.
Define $\mu\in{\rm Hom}_{\sZ}(\G_4,\F)$ as follows. For any
$\g\in\G_4$, we can write $\g=\tau+n\b$ with $n\in\Z$ and
$\tau\in\G_2$. Define $\mu_4(\g)=\mu''(\tau)+n\mu_3(\b)$. Suppose
$\tau+n\b=0$ for some $n\in\Z$. Then $\tau=-n\b\in\G_2\cap\G_3$. Since
$\mu''|_{\G_2\cap\G_3}=\mu_3|_{\G_2\cap\G_3}$, we have $\mu''(\tau)=
\mu_3(\tau)=\mu_3(-n\b)$.
But obviously, $\mu_3(-n\b)=-n\mu_3(\be)$. Hence $\mu''(\tau)+n\mu_3(\b)=0$.
This shows that $\mu_4\in{\rm Hom}_{\Z}(\G_4,\F)$ is uniquely defined.
So $d|_{\BB(\G_4,J)}=d_{u',\mu_4,f'_1,f'_2,f'_3,f'_4}$ and
$\G_4\supset\G_2$, $\G_4\ne\G_2$. This contradicts the maximality of $\G_2$.
Therefore, $\G_2=\G$ and $d=d_{u'',\mu'',f''_1,f''_2,f''_3,f''_4}$
satisfies the condition in (2.5).
\qed\par\ni
{\bf Lemma 2.5}. {\it Let $d$ be any derivation of $\BB$. Write
$d=\sum_{\a\in \G}d_\a$ with $d_\a\in(\mbox{\it Der}\:\BB)_\a.$
Then $d_\a=0$ for all but a finite $\a\in\G.$
}
\par\ni
{\it Proof}. By Lemmas 2.1-4, for any $\a\in\G\bs\{0,\si_1,\si_2\}$,
if $d_\a\ne0$,
there exists $u_\a\in\BB_\a$ with the leading degree, say, $\ii_\a$
such that $d_\a=\ad_{u_\a}$.
Fix $v=x^{\b,\jj}\in\BB$ with $\b_1\ne0$ or $j_1\ne0$. Then
$d(v)$ is contained in a sum of finite number of $\BB_\a$. Thus
$[x^{\a,\ii_\a},v]=0$ for all but a finite number of $\a$. But
by (2.8), we see that $[x^{\a,\ii},v]\ne0$ for all but a finite number of
$\a$. Thus $\{\a\,|\,d_\a\ne0\}$ is finite.
\qed\par\ni
{\bf Theorem 2.6}. {\it (1) $\der\BB=\oplus_{\a\in\G}(\der\BB)_\a,$
where each $(\der\BB)_\a$ consists of derivations of the form (2.40).
(2) The set $(\der\BB)_{\it f}$ of locally finite elements of $\der\BB$ is
\par\ni\hs{5pt}$
\left\{
\matrix{
\F\ad_1+{\it Hom}^*_{\sZ}(\G,\F)+\F\ptl_{t_2},\mbox{ \ or}
\vs{4pt}\hfill\cr
{\it span}_{\sF}\{\ad_{x^{\a,\ii}}\,|\,i_1\!=\!0\mbox{ or }
(\a,\ii)\!=\!(0,1_{[1]})\}\!+\!{\it Hom}_{\sZ}(\G,\F)
\!+\!\F d_1\!+\!\F\ptl_{t_2},
\hfill\cr}\right.
$\hfill(2.42)\par\ni
if $\pi_1(\G)\ne\{0\}$ or $\pi_1(\G)\!=\!\{0\}$ respectively.}
\par\ni
{\it Proof}. (1) follows from Lemma 2.5, and
clearly elements in (2.42) are locally finite.
\ul{\it Case (i)}: $\pi(\G)\ne\{0\}$.
Suppose $d=\sum_{\a\in\G_0}d_\a$ is locally finite, where $\G_0$ is a
finite subset of $\G$ such that all $d_\a\ne0,\a\in\G_0$.
Choose a
total ordering compatible with the group structure on $\G$.
Let $\b$ be the maximal element of $\G_0$. If $\b\ne0$ by reversing the
ordering if necessary, we can suppose $\b>0$, and also by (2.40),
we can suppose $d_\b=\ad_{u_\a}|_{\BB}$ for some $u_\a$.
If $\b_1\ne0$, then $d_\b=\ad_{u_\a}$ is inner, let $\ii$ be the leading
degree of $u_\a$ and say the coefficient of the leading term is $1$,
then we have
\par\ni\hs{5pt}$\dis
d^k(x^{2\b,0})\!\equiv\! \ad^k_{x^{\b,\ii}}(x^{2\b,0})
\!\equiv\! k!\b_1 x^{(k+2)\b,k\ii}\!\not\equiv\!0
\,({\rm mod}\sum_{\a<(k+2)\b}\BB_\a\!+\!
\BB_{(k+2)\b}^{(k\ii)}),
\hfill(2.43)$\par\ni
which shows that $d$ is not locally finite. Suppose $\b_1=0$.
Choose $\g\in\G$ with $\g_1\ne0$,
then
$$
d^k(x^{\g,0})\equiv \ad^k_{x^{\b,\ii}}(x^{\g,0})
\equiv (-(\b_2-1)\g_1)^k x^{k\b+\g,k\ii}\not\equiv0
\,({\rm mod\,}\sum_{\a<k\b+\g}\BB_\a+
\BB_{k\b+\g}^{(k\ii)}),
$$
which again shows that $d$ is not locally finite
if $\b_2\ne0$.
If $\b_2=1$, then $\ii\ne0$
(since $x^{\si_1,0}=0$ from the statement before (1.8)),
say $i_2>0$, then
\par\ni\hs{5pt}$
\matrix{\dis
d^k(x^{\g,0})\!\!\!\!&\equiv \ad^k_{x^{\b,\ii}}(x^{\g,0})
\vs{4pt}\hfill\cr&
\equiv (-i_2\g_1)^k x^{k\b+\g,k(\ii-1_{[2]}}\not\equiv0
\,({\rm mod\,}\sum_{\a<(k+2)\b}\BB_\a+
\BB_{k\b+\g}^{(k(\ii-1_{[2]}))}).
\hfill\cr}
\hfill(2.44)$\par\ni
Thus suppose $\G_0=\{0\}$. If $d$ is not in (2.42), write $d=d'+d''$, where
$d'$ is in (2.42) and $d''=\ad_{u_0}$ is inner with the leading degree of $u_0$
being $\ii>0$. Say again, $i_2>0$, then (2.44) again shows that
$d$ is not locally finite. This proves the theorem in this case.
\ul{\it Case (ii)}:  $\pi_1(\G)=\{0\}$.
We give a filtration on $\BB: \BB_{[0]}\subset\BB_{[1]}\subset...$
by defining
$$
\BB_{[m]}={\rm span}\{x^{\a,\ii}\,|\,(\a,\ii)\in\G^\#\times J
\mbox{ with } i_1\le m\}\ \for\ m\in\N.
\eqno(2.46)$$
Suppose $d$ is a derivation not in (2.42). Write $d=d'+d''$ with
$d'$ in (2.42) and $d''=\ad_u$ for some $u\in\BB$. Let $m$ be the smallest
number such that $u\in\BB_{[m]}$. Write $u=u_1+u_2$ such that
$u_1$ is in the space
${\rm span}\{x^{\a,\ii}\,|\,i_1=m\}$ and $u_2\in\BB_{[m-1]}$.
Them $m>1$ or $m=1$ but
$u_1\notin\F x^{0,1_{[1]}}$. Using arguments as before
for $u_1$ we can deduce that $d$ is not locally finite.
\qed
\par\
\vs{-5pt}\par \cl{\bf 3. \ PROOF OF THEOREM 1.2} Now we are ready
to prove Theorem 1.2.
\par\ni
{\it Proof of Theorem 1.2}.
As for the second statement, observe that ${\cal M}_1,...,{\cal M}_4$
correspond to the isomorphism classes of $\BB(\G,J)$ with
$J=\{0\},\N\times\{0\},\{0\}\times\N,\N^2$ respectively. So, we shall prove
the first statement below.
We shall ALWAYS use the same notations with a prime to denote elements
associated with $\BB'$.
\par
``$\Lar$'':
For convenience, we denote $t_1=x^{0,1_{[1]}},t_2=x^{0,1_{[2]}}$,
and $x^\a=x^{\a,0}$.
Define another algebra structure on $(\AA_2,\odot)$
(not necessarily associative) by:
$u\odot v=\ptl_1(u)(\ptl_2(v)-v).$
Then we have
\par\ni$
x^{\a,\ii}\odot x^{\b,\jj}=
\left\{\matrix{
x^{\a+\b,\ii+\jj-\sone}(\a_1t_1+i_1)((\b_2-1)t_2+j_2)
\mbox{ if }\ii+\jj-\one\in J,
\vs{4pt}\hfill\cr
\a_1x^{\a+\b,\ii+\jj-1_{[2]}}((\b_2-1)t_2+j_2)
\mbox{ if }i_1+j_1=0,i_2+j_2>0,
\vs{4pt}\hfill\cr
(\b_2-1)x^{\a+\b,\ii+\jj-1_{[1]}}(\a_1t_1+i_1)
\mbox{ if }i_1+j_1>0,i_2+j_2=0,
\vs{4pt}\hfill\cr
\a_1(\b_2-1)x^{\a+\b,\ii+\jj}
\mbox{ if }i_1+j_1=0,i_2+j_2=0,
\hfill\cr
}\right.
$\hfill(3.1)\par\ni
Then the Lie bracket $[\cdot,\cdot]$ is given by
$[u,v]=u\odot v-v\odot u.$
Obviously, using $\phi:(\b_1,\b_2)\mapsto
(a\b_1,\b_2+b\b_1)$, we see that $\si_i\in\G$ if and only if $\si_i\in\G'$
and $\phi(\si_i)=\si_i$ if it is in $\G$ for $i=1,2$.
\par
Suppose $b=0$. We shall verify that
$\psi: x^{\b,\jj}
\mapsto a^{-1}x'^{\b'}(at'_1)^{i_1}(t'_2)^{i_2},$
is an isomorphism $(\AA_2,\odot)\cong(\AA'_2,\odot)$,
where we use the notation $\b'=\phi(\b)$.
Say, we have the most complicated case, $\ii+\jj-\one\in J$.
Using (3.1), we have
$$
\matrix{
\psi(x^{\a,\ii}\odot x^{\b,\jj})
\!=\!a^{-1} x'^{\a'+\b'}
(at'_1)^{i_1+j_1-1}(t'_2)^{i_2+j_2-1}(\a_1at'_1+i_1)((\b_2-1)t'_2+j_2),
\vs{4pt}\hfill\cr
\psi(x^{\a,\ii})\odot\psi(x^{\b,\jj})
\!=\!a^{-2} x'^{\a'+\b'}
(at'_1)^{i_1+j_1-1}(t'_2)^{i_2\!+\!j_2\!-\!1}
(\a'_1at'_1\!+\!ai_1)((\b'_2\!-\!1)t'_2\!+\!j_2),
\hfill\cr}
$$
where in the second equation,
$a$ must appears before $i_1$ because it is arisen
from $\ptl'_1((at'_1)^{i_1})$.
It is immediately to check that
$(\a_1at'_1+i_1)((\b_2-1)t'_2+j_2)=
a^{-1}(\a'_1at'_1+ai_1)((\b'_2-1)t'_2+j_2),$
since $\a'_1=a\a_1,\b'_2=\b_2$.
Thus $\psi$ is an isomorphism
$(\AA_2,\odot)\cong(\AA'_2,\odot)$, which
induces an isomorphism of $\BB\cong\BB'$.
\par
Suppose $b\ne0$. We define
$\psi(x^{\b,\jj})=a^{-1}x'^{\b'}(at'_1+bt'_2)^{j_1}(t'_2)^{j_2}.$
We claim that it is an isomorphism $\BB\cong\BB'$ (but not necessarily an
isomorphism from $(\AA,\odot)$ to $(\AA',\odot)$), i.e., we want to
prove
$$
[\psi(x^{\a,\ii}),\psi(x^{\b,\jj})]=
\psi([x^{\a,\ii},x^{\b,\jj}]).
\eqno(3.2)$$
Again say, $\ii+\jj-\one\in J$.
First we calculate
$\psi(x^{\a,\ii}\odot x^{\b,\jj})$, which is the term
$$
a^{-1} x'^{\a'+\b'}
(at'_1+bt'_2)^{i_1+j_1-2}(t'_2)^{i_2+j_2-1},
\eqno(3.3)$$
where if $i_1+j_1-2<0$ the corresponding factor does not appear,
multiplied by the term
$$
(at'_1+bt'_2)(\a_1(at'_1+bt'_2)+i_1)((\b_2-1)t'_2+j_2).
\eqno(3.4)$$
Similarly,
$\psi(x^{\a,\ii})\odot\psi(x^{\b,\jj})$ is (3.3) multiplied by
$$
a^{-1}(\a'_1(at'_1+bt'_2)+ai_1)(((\b'_2-1)t'_2+j_2)(at'_1+bt'_2)+j_1bt'_2),
\eqno(3.5)$$
where the factor $(\a'_1(at'_1+bt'_2)+ai_1)$ is arisen from
$\ptl'_1(x'^{\a'}(at'_1+bt'_2)^{i_1})$ and the last factor is arisen
from $(\ptl'_2-1)(x'^{\b'}(at'_1+bt'_2)^{j_1}(t'_2)^{j_2})$.
Thus (3.2) is equivalent to
\def\sc#1{\rb{1pt}{\mbox{$\scriptstyle#1$}}}
$$
\matrix{
(at'_1\!+\!bt'_2)((\a_1{\sc(}at'_1\!+\!bt'_2{\sc)}\!+\!i_1)
({\sc(}\b_2\!-\!1{\sc)}t'_2\!+\!j_2)
\!-\!(\b_1{\sc(}at'_1\!+\!bt'_2{\sc)}\!+\!j_1)({\sc(}\a_2\!-\!1{\sc)}
t'_2\!+\!i_2))
\vs{4pt}\hfill\cr
=
a^{-1}((\a'_1(at'_1+bt'_2)+ai_1)
(((\b'_2-1)t'_2+j_2)(at'_1+bt'_2)+j_1bt'_2)
\vs{4pt}\hfill\cr\ \ \ \ \ \ \
-(\b'_1(at'_1+bt'_2)+aj_1)(((\a'_2-1)t'_2+i_2)(at'_1+bt'_2)+i_1bt'_2)),
\hfill\cr}
$$
\def\sc{\scriptstyle}
which is straightforward to verify.
\par
``$\Rar$'':
The isomorphism induces an isomorphism $\psi:\der\BB\rar\der\BB'$, which
maps locally nilpotent, semi-simple, locally finite elements to
locally nilpotent, semi-simple, locally finite elements respectively.
We consider in two cases.
\par
{\it Case 1}. $\pi_1(\G)\ne\{0\}$.
\par
By Theorem 2.6, $\BB$ has no {\it ad}-locally nilpotent element, and
$\F1$ are the only
{\it ad}-locally finite elements of $\BB$, thus
$\BB'$ has no {\it ad}-locally nilpotent element, and
$\F\psi(1)$ are the only
{\it ad}-locally finite elements of $\BB'$.
Hence $\pi_1(\G')\ne\{0\}$ and
$\psi(1)\in\F1'.$
Note that $\F1$ is the set of {\it ad}-semi-simple elements of
$\BB$ if and only if $J_1=\{0\}$. Thus $J_1=J'_1$.
Define $\ptl_{t_1}=\ad_1-\pi_1:x^{\a,\ii}\mapsto i_1x^{\a,\ii-1_{[1]}}$ if
$J_1\ne\{0\}$ and $\ptl_{t_1}=0$ if $J_1=\{0\}$.
Then by (2.42), we see that
$\F\ptl_{t_1}+\F\ptl_{t_2}$ is the set of locally nilpotent element of
$\der\BB$, since the dimension of
$\F\ptl_{t_1}+\F\ptl_{t_2}$ is the number of $p=1,2$ with $J_p\ne\{0\}$,
thus $J_2=J'_2$. This proves that $J=J'$.
We may assume that $J\ne\{0\}$ since the result follows from Refs.~(1) and (2)
if $J=\{0\}$. Thus $\G^\#=\G$ (cf.~statements before (1.8)).
\par
Since the nonzero vectors in
$\cup_{\a\in\G}\BB_\a$ are the only common eigenvectors for
the derivations in ${\rm Hom}_{\sZ}(\G,\F)$, which
is the set of semi-simple derivations of $\der\BB$ (cf.~(2.42)).
It follows that there exists a bijection $\phi:\G\rar\G'$ such that
$\psi(\BB_\a)=\BB'_{\phi(\a)}$ for $\a\in\G,$ and $\phi(0)=0.$
Thus also
$\psi((\der\BB)_\a)=(\der\BB')_{\phi(\a)}$ for $\a\in\G$ (cf.~(2.4)).
Since $0\ne[(\der\BB)_\a,(\der\BB)_\b]\subset (\der\BB)_{\a+\b}$ for all
$\a,\b\in\G$, we obtain that
$\phi(\a+\b)=\phi(\a)+\phi(\b)$, i.e., $\phi$
is a group isomorphism.
Let $\a\in\G\bs\{\si_1\}$ be such that $\a'=\phi(\a)\in\G'\bs\{\si_1\}$.
Then since $\F x^{\a,0}=\{u\in\BB_\a\,|\,(\F\ptl_1+\F\ptl_2)(u)=0\}$,
we obtain that
$\psi(x^{\a,0})=a_\a x'^{\a',0}$ for some $a_\a\in\F\bs\{0\}.$
Applying $\psi$ to $[x^{\a,0},x^{\b,0}]=
(\a_1(\b_2-1)-\b_1(\a_2-1))x^{\a+\b,0}$, we obtain
$$
a_\a a_\b(\a'_1(\b'_2-1)-\b'_1(\a'_2-1))=
a_{\b+\a}(\a_1(\b_2-1)-\b_1(\a_2-1)),
\eqno(3.6)$$
for $\a,\b,\a+\b,\a',\b',\a'+\b'\ne\si_1$.
Take $\b=0$, this gives
$$
\a'_1=a\a_1 \mbox{ \ for all }\ \a\in\G,
\mbox{ \ where \ }a=a_0^{-1}.
\eqno(3.7)$$
Assume that $\si_1\in\G$.
If $J_1=\{0\}$ or $J_2=\{0\}$, then $\ol d_1$ or $d_1$ is an outer
derivation in $(\der\BB)_{\si_1}$
and $(\der\BB)_\a=\ad_{\BB_\a}$ for all $\a\in\G\bs\{0,\si_1\}$,
thus $\si_1$ must be in $\G'$ and $\phi(\si_1)=\si_1$.
If $J_1\ne\{0\}\ne J_2$, then
$\si_1$ is the only $\a$ such that
$\{u\in\BB_\a\,|\,(\F\ptl_{t_1}+\F\ptl_{t_2})(u)=0\}$ has dimension
2 (spanned by $x^{\si_1,1_{[p]}},p=1,2$).
We must have $\si_1\in\G'$ and $\phi(\si_1)=\si_1$.
This proves that $\si_1\in\G\Rla\si_1=\phi(\si_1)\in\G'$.
Suppose $\si_2\in\G'$ but $\g=\phi^{-1}(\si_2)\ne\si_2$. By (3.7),
$\g_1=0$, and so $\g_2\ne2$ since $\g\ne\si_2$ (cf.~(1.6)).
Take $\a\in\G\bs\{0\}$ with $\a_1\ne0$, then by applying $\psi$ to
$[x^{\a,0},x^{\g-\a,0}]=\a_1(\g_2-2)x^{\g,0}\ne0,$
we obtain that $[x^{\a',0},x^{\si_2-\a',0}]\ne0$, but by (2.8)
it is zero, a contradiction. Thus
$\phi(\si_2)=\si_2$.
This proves that $\si_2\in\G\Rla\si_2=\phi(\si_2)\in\G'$.
In (3.6), let $\b=m\a,m\ge2,\a_1\ne0$, then we obtain
$a_0a_{(m+1)\a}=a_\a a_{m\a}$, from this we obtain
$$
a^{m-1}_0a_{m\a}=a_\a^m\mbox{ for }m\in\Z.
\eqno(3.8)$$
We define $a_{\si_1}=a_0^2a_{-\si_1}^{-1}$ if $\si_1\in\G$ and
$a_{\si_2}=a_0^2a_{-\si_2}^{-1}$ if $\si_2\in\G$. Using
(3.8), we can prove that (3.6) holds for all $\a,\b\in\G$.
For example, if $\a=\si_1$, then (3.6) holds trivially. Suppose
$\a,\b\ne\si_1$ and $\a+\b=\si_1$, then
$\b_1=-\a_1,\b_2=-\a_2+1$ and (3.6) is equivalent to
$a_\a a_\b=a_0a_{\si_1}$ which is equivalent to
$a_0^4a^{-1}_{-\a}a^{-1}_{-\b}=a_0a_0^2a_{-\si_1}^{-1}$ which indeed
holds by setting $\a$, $\b$ to be $-\a,-\b$ in (3.6).
Now as proved in Ref.~(2),
or as an exercise by using (3.7) to prove it, we have
$a_\a a_\b=a_0 a_{\a+\b}$ for $\a,\b\in\G.$
Fix any $\a\in\G$ with $\a_1\ne0$ and set $b=
\a_2\a_1^{-1}(\a'_2\a_2^{-1}-1)$, then this, together with (3.6), (3.7),
shows that $\b'_2=\b_2+b\b_1$. This proves that $\phi:(\b_1,\b_2)\mapsto
(a\b_1,\b_2+b\b_1)$ has the required form.
\par
Suppose $J_1\ne\{0\}=J_2$. Then we have
$\psi(x^{\a,1_{[1]}})=b_\a x'^{\a',1_{[1]}}+c_\a x'^{\a',0}$
for some $b_\a,c_\a\in\F$
and all $\a\in\G\bs\{\si_1\}.$
This follows from that $\F x^{\a,1_{[1]}}+\F x^{\a,0}$ is the
subspace $\{u\in\BB_\a\,|\,\ptl_{t_1}^2(u)=0\},$
and that $\F\ptl_{t_1}$
are the only locally nilpotent derivations of $\BB$. Applying $\psi$ to
$[x^{\a,1_{[1]}},x^{\b,0}]=(\a_1(\b_2-1)-\b_1(\a_2-1))x^{\a+\b,1_{[1]}}
+(\b_2-1)x^{\a+\b,0},$
we obtain
$$
\matrix{
b_\a a_\b(\a'_1(\b'_2-1)-\b'_1(\a'_2-1))
=b_{\a+\b}(\a_1(\b_2-1)-\b_1(\a_2-1))
\vs{4pt}\hfill\cr
b_\a a_\b(\b'_2-1)
+c_\a a_\b(\a'_1(\b'_2-1)-\b'_1(\a'_2-1))
\hfill\cr
\ \ \ \ \ \ \ \ \
=a_{\a+\b}(\b_2-1)+c_{\a+\b}(\a_1(\b_2-1)-\b_1(\a_2-1)),
\hfill\cr}
\eqno\matrix{(3.9)\vs{4pt}\cr\cr(3.10)\cr}$$
for all $\a,\b,\a+\b\in\G\bs\{\si_1\}.$
Let $\a=0$ in (3.9) and let $\b=\a$ in (3.10), we obtain
$$
b_0 a_\b\b'_1=b_\b\b_1,\ \
b_\a a_\a(\a'_2-1)=a_{2\a}(\a_2-1).
\eqno(3.11)$$
Thus we have $b_0=1$ by letting $\a=0$ in the last equation of (3.11)
and hence by the first of (3.11),
$b_\a=\a'_1\a_1^{-1}a_\a=a^{-1}_0a_\a$ by (3.7)
and $b_\a a_\a=a^{-1}_0a_\a^2=a_{2\a}$ by (3.8).
Thus the last equation of (3.11) gives $\a'_2=\a_2$ for all $\a\in\G$.
Thus $b=0$ because $\phi:\a\mapsto\a'=(a\a_1,\a_2+b\a_1)$.
This proves the theorem in this case.
\par
{\it Case 2}. $\pi_1(\G)=\{0\}$.
\par
By Case 1, $\pi_1(\G')=\{0\}$. So $J_1=J'_1=\N$.
Using notation (2.46), by (2.42), the set
of {\it ad}-locally nilpotent elements of $\BB$
is $\BB_{[0]}$, thus $\psi(\BB_{[0]})=\BB'_{[0]}.$
Note that the common eigenvectors
for the derivations in (2.42) are the nonzero elements of the set
$$
\left\{\matrix{
\dis \bigcup_{\a\in\G}\F x^{\a,0},
&\mbox{if \ }\si_1\notin\G,
\vs{4pt}\hfill\cr
\dis\bigcup_{\a\in\G\bs\{\si_1\}}\F x^{\a,0}
\bigcup {\rm span}\{x^{\si_1,\ii}\,|\,|\ii|\le1\},
&\mbox{if \ }\si_1\in\G.
\hfill\cr}\right.
\eqno\matrix{^{\dis(3.12)}\vs{8pt}\cr(3.13)\cr}$$
\hs{3ex}
First assume that $\si_1\in\G$ but $\si_1\notin\G'$. Then we must have
$J_2=\{0\}$, otherwise if $J_2=\N$,
we would have that ${\rm span}\{x^{\si_1,\ii}\,|\,
|\ii|\le1\}$ in (3.13) were two dimensional but all
$\F x'^{\a,0}$ in $\BB'$ corresponding to (3.12) are one dimensional;
this is a contradiction.
Furthermore we have $J'_2\ne\{0\}$, otherwise
if $J'_2=\{0\}$, then $d'_1=\ptl'_{t_2}=0$ and by (2.42),
the set $(\der\BB')_{\rm f}$ of locally finite elements of $\BB'$
would be spanned by
inner locally finite derivations and semi-simple derivations, but
this is not true for $(\der\BB)_{\rm f}$, again we get a contradiction.
Now we have an outer locally nilpotent
derivation $d_1$ of $\BB$ which kills the set $\BB_{[0]}$ (since by (2.12),
$d_1(x^{\b,\jj})=-j_1 x^{\b,\jj-1_{[1]}}$ in this case),
but we do not
have an outer locally nilpotent derivation of $\BB'$ which kills $\BB'_{[0]}$
(since $J'_2\ne\{0\}$, $\ptl'_{t_2}$ does not kill $\BB'_{[0]}$), again
a contradiction. This proves that $\si_1\in\G\Rla\si_1\in\G'$.
\vs{-1pt}\par
Thus by (3.12) or (3.13), there exists a bijection $\phi:\a\rar\a'$
from $\G\rar\G'$ such
\vs{-1pt}that
$$
\psi(x^{\a,0})=a_\a x'^{\a',0}\mbox{ for some }a_\a\in\F\bs\{0\}
\mbox{ and }\a\in\G\bs\{\si_1\}. \vs{-1pt}\eqno(3.14)$$ If
$J_2=\{0\}\ne J'_2$, then $\BB_{[0]}$ is spanned by (3.12) or
(3.13), but $\BB'_{[0]}$ which contains $x'^{0,2_{[2]}}$, is not
spanned by the corresponding set of (3.12) or (3.13) in $\BB'$.
Thus we obtain that $J_2=\{0\}\Rla J'_2=\{0\}$. This proves that
$J=J'$. \vs{-1pt}\par Applying the action of (2.42) to (3.12) or
(3.13), if we denote by $N_\a$ the elements of (2.42) which kill
$x^{\a}$ if $\a\ne\si_1$, or which kill
$\{x^{\si_1,\ii}\,|\,|\ii|\le1\}$ if $\a=\si_1$, then $N_\a$ has
codimension 1 in (2.42), mainly $x^{0,1_{[1]}}\notin N_\a$, if
$\a=0$; or has codimension at least 2, mainly $x^{0,1_{[1]}}\notin
N_\a$ and $d_\mu\notin N_\a$ for all $\mu$ satisfying
$\mu(\a)\ne0$, if $\a\ne0$. Thus we can suppose $\psi(1)=b_0 1'$
for some $b_0\in\F\bs\{0\}$. Observe from (2.42) that all {\it
ad}-semi-simple elements must take the form $a_0
x^{0,1_{[1]}}+\sum_{(\a,\ii)\in G_0} c_{\a,\ii}x^{\a,\ii},$ where
$a_0\ne0$ and $G_0$ is a finite subset of $\G\times J$ such that
$i_1=0$ if $(\a,\ii)\in\G_0$ or $|\ii|\le1$ if $\a=\si_1$. Thus,
we can suppose $\psi(x^{0,1_{[1]}})=a'_0 x'^{0,1_{[1]}}
+\sum_{(\a',\ii)\in G'_0} c_{\a',\ii}x'^{\a',\ii}.$ Applying
$\psi$ to $[1,x^{0,1_{[1]}}]=1$ we obtain $[b_01',a'_0
x'^{0,1_{[1]}} +\sum_{(\a',\ii)\in G'_0}
c_{\a',\ii}x'^{\a',\ii}]=b_01',$ thus $a'_0=1$. Finally applying
$\psi$ to $[x^{0,1_{[1]}},x^{\b,0}]=(\b_2-1)x^{\b,0}$ using
(3.14), we obtain $[x'^{0,1_{[1]}}+\sum_{(\a',\ii)\in G'_0}
c_{\a',\ii}x'^{\a',\ii}, a_\b x'^{\b',0}]=(\b_2-1)a_\b
x'^{\b',0}.$ This shows that $(\b'_2-1)=\b_2-1$, i.e.,
$\b'_2=\b_2$, and so, $\G=\G'$ and $\phi$ is the identity
isomorphism. This completes the proof of the theorem. \qed
\small\par\ \vs{-5pt}\par \cl{\bf ACKNOWLEDGEMENTS} \vs{-1pt}\par
Part of this work was conducted during the first author's visit to
The Hong Kong University of Science and Technology. He would like
to thank Dr.~Xiaoping Xu for hospitality and help. This work is
supported by a NSF Fund of China and a Fund from National
Education Ministry of China.
\par\
\vs{-5pt}\par \cl{\bf REFERENCES}
\par\ni\hi3ex\ha1
1.~Block,~R. On torsion-free abelian groups and Lie algebras.
  Proc.~Amer. Math.~Soc. {\bf1958}, {\it9}, 613-620.
  \par\ni\hi3ex\ha1
2.~Dokovic,~D.; Zhao,~K. Derivations, isomorphisms and
  second cohomology of generalized Block algebras.
  Algebra Colloquium {\bf1996}, {\it3}, 245-272.
  \par\ni\hi3ex\ha1
3.~Dokovic,~D.; Zhao,~K. Some infinite dimensional simple Lie
  algebras related to those of Block.
  J.~Pure Appl.~Alg. {\bf1998}, {\it127}, 153-165.
  \par\ni\hi3ex\ha1
4.~Dokovic,~D.; Zhao,~K. Some simple subalgebras of generalized
  Block algebras.
  J.~Alg. {\bf1997}, {\it192}, 74-101.
  \par\ni\hi3ex\ha1
5.~Mazirchuk,~V. Verma modules over generalized Witt algebra.
  {\it Compositio Mathematica} {\bf1999}, {\it115}, 21-35.
  \par\ni\hi3ex\ha1
6.~Osborn,~J.~M.; Zhao,~K. Infinite-dimensional Lie algebras
  of generalized Block type.
  Proc. Amer.~Math.~Soc. {\bf1999}, {\it127}, 1641-1650.
  \par\ni\hi3ex\ha1
7.~Su,~Y.; Xu,~X. Structure of divergence-free Lie algebras.
  J.~Alg., {\bf2001}, {\it243}, 557-595.
  \par\ni\hi3ex\ha1
8.~Su,~Y.; Xu,~X.; Zhang,~H. Derivation-simple algebras
  and the structures of Lie algebras of Witt type.
  J.~Alg. {\bf2000}, {\it233}, 642-662.
  \par\ni\hi3ex\ha1
9.~Xu,~X. New generalized simple Lie algebras of Cartan type over a
  field with characteristic 0.
  J.~Alg. {\bf2000}, {\it224}, 23-58.
  \par\ni\hi3.2ex\ha1
10.~Xu,~X. On simple Novikov algebras and their irreducible modules.
  J.~Alg. {\bf1996}, {\it185}, 905-934.
  \par\ni\hi3.2ex\ha1
11.~Xu,~X. Novikov-Poisson algebras.
  J.~Alg. {\bf1997}, {\it190}, 253-279.
  \par\ni\hi3.2ex\ha1
12.~Xu,~X. Generalizations of Block algebras.
   Manuscripta Math. {\bf1999}, {\it100}, 489-518.
  \par\ni\hi3.2ex\ha1
13.~Xu,~X. Quadratic conformal superalgebras.
   J.~Alg. {\bf2000}, {\it224}, 1-38.
  \par\ni\hi3.2ex\ha1
14.~Zhao,~K. A Class of infinite dimensional simple Lie algebras.
  J.~London Math.~Soc. (2), {\bf2000}, {\it62}, 71-84.
\end{document}